\documentclass[11pt]{amsart}
\usepackage{amscd,amsfonts,amsmath,amssymb,amsthm}
\usepackage{booktabs,cases,cite,color,diagbox}
\usepackage{enumitem,fancyhdr,graphicx}
\usepackage{accents,hyperref,indentfirst}
\usepackage{leftidx,mathrsfs,multirow}
\usepackage{picins,ragged2e,rotating,slashbox}
\usepackage{verbatim}
\usepackage{float}
\usepackage[margin=1.2in]{geometry}
\usepackage[all]{xy}
\usepackage{bm}
\newtheorem{definition}{Definition}
\newtheorem{theorem}{Theorem}
\newtheorem{lemma}{Lemma}
\newtheorem{proposition}{Proposition}
\newtheorem{corollary}{Corollary}
\newtheorem{theoremA}{Theorem}

\begin{document}

\title[The Orlicz-Gauss image problem for pseudo-cones]{The Orlicz-Gauss image problem for pseudo-cones and its associated spherical optimal transport}

\date{}

\author{Siqi Lei \& Xudong Wang}

\begingroup
\renewcommand{\thefootnote}{}
\footnotetext{School of Mathematics and Statistics, Shaanxi Normal University, Xi'an, 710119, China.}
\footnotetext{Email: sqlei@snnu.edu.cn (S. Lei), xdwang@snnu.edu.cn (X. Wang).}
\endgroup

\maketitle

\begin{abstract}
Pseudo-cones serve as the noncompact counterpart of convex bodies in convex geometry. This paper establishes a necessary and sufficient condition for the existence of solutions to the Orlicz-Gauss image problem for pseudo-cones and further demonstrates its connection to spherical optimal transport. Our approach combines the variational method with a novel restrictive technique, thereby strengthening the original result of Schneider up to a constant factor.
\end{abstract}

\section{Introduction}

Let $K\subset\mathbb{R}^n$ be a compact convex set with the origin in its interior, and $\lambda$ be a nonzero finite measure defined on spherical Lebesgue measurable sets that vanishes on sets of Hausdorff dimension $(n-2)$. The Orlicz-Gauss image measure of $K$ with respect to $\lambda$ is a spherical Borel measure given by the push-forward of the measure $\psi(\rho_K)\lambda$ via the reverse radial Gauss map $\alpha^*_K$, i.e.,
\begin{equation*}
\lambda_\psi(K,\cdot)=(\alpha^*_K)_\sharp\,\psi(\rho_K)\lambda,
\end{equation*}
where $\rho_K$ is the radial function of $K$ and $\psi$ is a continuous function on $(0,+\infty)$. When $\lambda$ is the spherical Lebesgue measure and $\psi\equiv1$, the Orlicz-Gauss image measure of $K$ is also known as the integral Gauss curvature of $K$, and written as $J(K,\cdot)$, which is one of the most fundamental geometric measures in convex geometry.

\subsection{The Aleksandrov problem}

The classical Aleksandrov problem \cite{Aleksandrov-Existence_and_uniqueness_of_a_convex_surface_with_a_given_integral_curvature} is to find a compact convex set with prescribed integral Gauss curvature. That is, for a Borel measure $\mu$ on $\mathbb{S}^{n-1}$, what are the necessary and sufficient conditions on $\mu$ in order that there exists a compact convex set $K$ with the origin in its interior such that $J(K,\cdot)=\mu$? It is characterized by the famous {\it Aleksandrov's condition}:
\begin{equation*}
\begin{cases}
\mu(\mathbb{S}^{n-1})=|\mathbb{S}^{n-1}|,\\
\mu(\omega)<|\mathbb{S}^{n-1}\setminus\omega^*|,
\end{cases}
\end{equation*}
where $\omega^*=\{v\in\mathbb{S}^{n-1}:\langle u,v\rangle\leqslant0,\ \mbox{for all}\ u\in\omega\}$ is the polar set of the spherically convex set $\omega\subset\mathbb{S}^{n-1}$.

Oliker \cite{Oliker-Embedding_S_n_into_R_n_1_with_given_integral_Gauss_curvature} demonstrated that the Aleksandrov problem belongs to a class of spherical optimal transport problems with the cost function
\begin{equation*}
c(u,v)=
\begin{cases}
\log\langle u,v\rangle, & \text{if}\ \,\langle u,v\rangle>0,\\
-\infty, & \text{otherwise}.
\end{cases}
\end{equation*}
For more details regarding the Aleksandrov problem and spherical optimal transport, we refer to \cite{Bertrand-Prescription_of_Gauss_curvature_using_optimal_mass_transport,Bertrand-On_the_Gauss_image_problem} and the references therein.

It is worth noting that a decomposition of $J(K,\cdot)$ was established in \cite{Wang-Zhu-The_decompositions_of_curvature_measures}. In particular, the density of $J(K,\cdot)$ with respect to the spherical Lebesgue measure is
\begin{equation*}
\frac{\rho^n_K(u)H_{n-1}(K,u)}{h_K(\alpha_K(u))}\ \ \text{for a.e.}\ \ u\in\mathbb{S}^{n-1},
\end{equation*}
where $h_K$ is the support function of $K$, $\alpha_K$ is the radial Gauss map of $K$, and $H_{n-1}(K,u)$ is the generalized Gauss curvature of $K$ at $\rho_K(u)u$. Therefore, the Monge-Amp\`{e}re type equation corresponding to the Aleksandrov problem is
\begin{equation*}
h\det(\nabla^2h+h\mathrm{I})=f(h^2+|\nabla h|^2)^\frac{n}{2},
\end{equation*}
where $\nabla$ is the Levi-Civita connection of the standard round metric $\mathrm{I}$ on $\mathbb{S}^n$. For the regularity of the Aleksandrov problem, please see, e.g., Guan and Li \cite{Guan-Li-C_1_1_estimates_for_solutions_of_problem_Alexandrov} and Li, Sheng, and Wang \cite{Li-Sheng-Wang-Flow_by_Gauss_curvature_to_the_Aleksandrov_and_dual_Minkowski_problems}. For the general Aleksandrov type problem, one can refer to \cite{Huang-LYZ-L_p_Aleksandrov_problem,Haberl-LYZ-The_even_Orlicz_Minkowski_problem,
Huang-LYZ-Geometric_measures_in_the_dual_Brunn_Minkowski_theory,Boroczky-Lutwak-Yang-Zhang-Zhao-The_Gauss_image_problem,
Li-Sheng-Wang-Flow_by_Gauss_curvature_to_the_Aleksandrov_and_dual_Minkowski_problems,
Feng-He-The_Orlicz_Aleksandrov_problem_for_Orlicz_integral_curvature,Gardner-General_volumes_in_the_Orlicz_Brunn_Minkowski_theory,
Huang-Xing-Ye-Zhu-On_the_Musielak_Orlicz_Gauss_image_problem} and the references therein.

\subsection{Pseudo-cones}

The concept of pseudo-cones stems from the works in \cite{Artstein-Avidan-A_zoo_of_dualities,Schneider-Pseudo_cones,Xu-Li-Leng-pseudo-cones}. Let $K\subset\mathbb{R}^n$ be a nonempty closed convex set not containing the origin, then $K$ is called a pseudo-cone if $\lambda K\subset K$ for any $\lambda>1$. The recession cone of $K$ is defined by
\begin{equation*}
\operatorname{rec}K=\{z\in\mathbb{R}^n:K+z\subset K\}.
\end{equation*}
If $\operatorname{rec}K$, denoted by $C$, is $n$-dimensional and pointed, then we call $K$ a $C$-pseudo-cone. The class of all $C$-pseudo-cones is denoted by $ps(C)$.

Throughout the paper, $C$ represents an $n$-dimensional and pointed closed convex cone. Moreover, we denote
\[
\Omega_C=\mathbb{S}^{n-1}\cap\operatorname{int}C\ \ \text{and}\ \ \Omega_{C^\circ}=\mathbb{S}^{n-1}\cap\operatorname{int}C^\circ,
\]
where the polar cone $C^\circ$ of $C$ is defined by $C^\circ=\{x\in\mathbb{R}^n:\langle x,y\rangle\leqslant0,\forall\,y\in C\}$.

For $K\in ps(C)$, the radial function of $K$ is defined by
\begin{equation*}
\rho_K(u)=\min\{r>0:ru\in K\},\ u\in \Omega_{C,K},
\end{equation*}
where $\Omega_{C,K}=\{u\in\operatorname{cl}\Omega_C:\text{there exists a}\ r>0\ \text{such that}\ ru\in K\}$. The radial Gauss image $\bm{\alpha}_K(\eta)$ of a Borel set $\eta\subset\Omega_C$ is the set of all unit outer normals of $K$ at $\rho_K(u)u$ for $u\in\eta$. For a vector $v\in\Omega_{C^\circ}$, the reverse Gauss map $\alpha^*_K(v)$ is the unique point $u\in\Omega_C$ such that $v$ is a unit outer normal of $K$ at $\rho_K(u)u$.

Let $\lambda$ be a nonzero finite measure on spherical Lebesgue measurable subsets of $\operatorname{cl}\Omega_{C^\circ}$, which is zero on sets of Hausdorff dimension $n-2$. For a positive continuous function $\phi$ on $(0,+\infty)$, the Orlicz-Gauss image measure of $K\in ps(C)$ is defined by $\lambda_\phi(K,\cdot)=(\alpha^*_K)_\sharp\,\phi(\rho_K)\lambda$, i.e.,
\begin{equation*}
\lambda_\phi(K,\eta)=\int_{\bm{\alpha}_K(\eta)}\phi(\rho_K(\alpha^*_K(v)))\,\mathrm{d}\lambda(v),\,\ \mbox{for any Borel set}\,\ \eta\subset\Omega_C.
\end{equation*}
Here, $\alpha^*_K(v)$ is well-defined for a.e. $v\in\bm{\alpha}_K(\Omega_C)$. If $\phi\equiv1$, then $\lambda_\phi(K,\cdot)=\lambda(K,\cdot)$ is called the Gauss image measure of $K$. In addition, if $\lambda$ is the spherical Lebesgue measure, then $\lambda(K,\cdot)=J(K,\cdot)$ is called the integral Gauss curvature of $K$.

The Orlicz-Gauss image problem for pseudo-cones reads: {\it Given a nonzero Borel measure $\mu$ on $\Omega_C$, what are the necessary and sufficient conditions on $\lambda$ and $\mu$ in order that there exists a $K\in ps(C)$ such that}
\begin{equation*}
\lambda_\phi(K,\cdot)=\mu?
\end{equation*}
For $\phi\equiv1$, if $\lambda$ is the spherical Lebesgue measure and $\mu$ is finite and concentrated on a compact domain of $\Omega_C$, Li, Ye, and Zhu \cite{Li-Zhu-Ye-The_dual_Minkowski_problem_for_unbounded_closed_convex_sets} proved that there exists $K\in ps(C)$ such that $J(K,\cdot)=\mu$ by the variational scheme. Generally, the Gauss image problem for $\lambda(K,\cdot)$ was investigated and completely characterized by Schneider \cite{Schneider-The_Gauss_image_problem_for_pseudo_cones}.

For further Minkowski-type problems pertaining to pseudo-cones, one can refer to e.g., \cite{Schneider-A_weighted_Minkowski_theorem_for_pseudo_cones,Schneider-The_Gauss_image_problem_for_pseudo_cones,
Li-Zhu-Ye-The_dual_Minkowski_problem_for_unbounded_closed_convex_sets,Schneider-Minkowski_type_theorems,
Schneider-A_Brunn_Minkowski_theory_for_coconvex_sets,Semenov-Zhao,Wang-Xu-Zhou-Zhu-Asymptotic_theory,Yang-Ye-Zhu-On_the_Lp}.

\subsection{Measure transport via pseudo-cones}

Throughout the paper, $\lambda$ is a nonzero finite measure on the Lebesgue measurable subsets of $\operatorname{cl}\Omega_{C^\circ}$, which is zero on sets of Hausdorff dimension $n-2$. We state Schneider's result on the Gauss image problem for pseudo-cones as follows.
\begin{theoremA}[see \cite{Schneider-The_Gauss_image_problem_for_pseudo_cones}]\label{Schneider-theorem}
Let $\mu$ be a nonzero Borel measure on $\Omega_C$, then there exists a $K\in ps(C)$ such that $\lambda(K,\cdot)=\mu$ if and only if $\mu(\Omega_C)=\lambda(\Omega_{C^\circ})$.
\end{theoremA}

Note that $\lambda(K,\cdot)=\mu$ also means $(\alpha^*_K)_\sharp\,\lambda=\mu$. Thus, if $\mu$ and $\lambda$ are probability measures on $\Omega_C$ and ${\rm cl}\,\Omega_{C^\circ}$ respectively, then there exists a $K\in ps(C)$ such that $(\alpha^*_K)_\sharp\,\lambda=\mu$. Exhilarating, Theorem \ref{Schneider-theorem} bears a strong resemblance to the celebrated Brenier-McCann theorem.
\begin{theoremA}[see \cite{Brenier-Polar_factorization_and_monotone_rearranoement_of_vector_valued_functions,
McCann-Existence_and_uniqueness_of_monotone_measure_preserving_maps}]\label{Brenier-McCann-theorem}
Let $\widetilde{\lambda},\widetilde{\mu}$ be two probability measure on $\mathbb{R}^n$, and $\widetilde{\lambda}$ be zero on sets of Hausdorff dimension $n-1$. There exists a convex function $f:\mathbb{R}^n\to(-\infty,+\infty]$ such that
\[
(\nabla f)_\sharp\widetilde{\lambda}=\widetilde{\mu},
\]
where $\nabla f$ is $\widetilde{\lambda}$-a.e. defined on $\operatorname{dom}f=\{x\in\mathbb{R}^n:f(x)<+\infty\}$.
\end{theoremA}

As Schneider \cite{Schneider-Pseudo_cones_and_measure_transport} pointed out, the gradient map $\nabla f$ in Theorem \ref{Brenier-McCann-theorem} corresponds to $\alpha^*_K$ in Theorem \ref{Schneider-theorem}. This shows that the Gauss image problem for pseudo-cones is actually a spherical optimal transport problem. Choosing a cost function $c:\Omega_{C^\circ}\times\Omega_C\to(0,+\infty)$ by
\begin{equation}\label{The-cost-function}
c(v,u)=\log\frac{1}{|\langle v,u\rangle|},
\end{equation}
Schneider proved that $\alpha^*_K$ is a maximizer of the supremum
\[
\sup_{T\in\mathcal{T}}\int_{\Omega_{C^\circ}}c(v,T(v))\,\mathrm{d}\lambda(v),
\]
where $\mathcal{T}$ consists of all measurable mappings $T$ from $\Omega_{C^\circ}$ to $\Omega_C$, defined $\lambda$-a.e., such that $T_\sharp\lambda=\mu$.

Further details on pseudo-cones and spherical optimal transport can be found in \cite{Schneider-Pseudo_cones_and_measure_transport,Schneider-Measure_transport_via_pseudo_cones,Schneider-Radial_measures_of_pseudo_cones}. For a more systematic treatment of spherical optimal transport, we recommend the recent important work by Li and Wu \cite{Li-Wu-Optimal_transports_on_the_spherical_domains_with_boundary_and_their_applications}.

\subsection{Our main results}

We develop a simple yet novel technique (termed the restrictive technique) which is different from Schneider's approximation approach. Consequently, we strengthen Theorem \ref{Schneider-theorem} up to a constant factor, extend it to the general Orlicz space, and provide the corresponding optimal transport interpretation. We again emphasize that this restrictive technique does not rely on uniform estimates, while Schneider's approximation approach does.

Throughout the paper, $\phi$ is a positive continuous function on $(0,+\infty)$. Let $K\in ps(C)$, $K$ is called a $C$-full set if $C\setminus K$ is bounded. In fact, the research of $C$-full sets originated from algebraic geometry \cite{Khovanskii-Timorin-On_the_theory_of_coconvex_bodies}. Using the variational method and the restrictive technique, we characterize the Orlicz-Gauss image problem as follows.
\begin{theorem}\label{The-Main-Theorem}
Let $\mu$ be a nonzero Borel measure on $\Omega_C$, then there exists a $C$-full set $K$ and a constant $\beta>0$ such that
\begin{equation*}
\beta\lambda_\phi(K,\cdot)=\mu
\end{equation*}
if and only if $\mu$ is finite. Moreover, such $\beta$ and $K$ are not unique. In fact, there exist infinitely many such pairs.
\end{theorem}

As by-products, we obtain the following corollaries.
\begin{corollary}[Enhanced Theorem \ref{Schneider-theorem} up to a constant]\label{The-corollary-1}
Let $\mu$ and $\lambda$ are probability measures on $\Omega_C$ and ${\rm cl}\,\Omega_{C^\circ}$ respectively, then there exists a $C$-full set $K$ and a constant $\beta\in(0,1]$ such that
\[
(\alpha^*_K)_\sharp\,\lambda=\beta\mu.
\]
\end{corollary}

When $\mu$ is finite and $\phi(t)=t^p$, Li and Li \cite{Li-Li-The-L_p_Gauss_image_problem_for_pseudo_cones} studied the $L_p$-Gauss image problem for $\lambda_p(K,\cdot)=\lambda_\phi(K,\cdot)$. Naturally, we cover and enhance their results as follows.
\begin{corollary}\label{The-corollary-2}
Let $\mu$ be a nonzero Borel measure on $\Omega_C$, $0\neq p\in\mathbb{R}$, then there exists a $C$-full set $K$ such that $\lambda_p(K,\cdot)=\mu$ if and only if $\mu$ is finite.
\end{corollary}

Let $K\in ps(C)$ and $\mu$ be a nonzero Borel measure on $\Omega_C$, we call $T$ a $(K,\phi)$-transport map from $\Omega_{C^\circ}$ to $\Omega_C$ if $\mu$ is the push-forward of $\phi(\rho_K)\lambda$, i.e., $T_\sharp\phi(\rho_K)\lambda=\mu$.
\begin{corollary}\label{The-corollary-3}
Let $\mu$ be a nonzero finite Borel measure on $\Omega_C$, then there exists a $C$-full set $K$ and a constant $\beta>0$ such that $\alpha^*_K$ is a maximizer of the supremum
\[
\sup_{T\in\mathcal{T}}\int_{\Omega_{C^\circ}}c(v,T(v))\,\mathrm{d}\lambda(v).
\]
Here $\mathcal{T}$ consists of all $(K,\widetilde{\phi})$-transport map from $\Omega_{C^\circ}$ to $\Omega_C$, $\widetilde{\phi}=\beta\phi$, and the cost function $c$ is as \eqref{The-cost-function}.
\end{corollary}

For the case of infinite measures, there are many difficulties in the Minkowski-type problem for pseudo-cones. Please see \cite{Schneider-Pseudo_cones} and \cite{Semenov-Zhao} for the in-depth study. Here, we obtain the following result for the infinite measures. Recall that a measure $\mu$ on $\Omega_C$ is called locally finite if $\mu(\eta)<+\infty$ for any compact $\eta\subset\Omega_C$.
\begin{theorem}\label{The-Main-Theorem-local-finite}
Let $\mu$ be a nonzero and locally finite Borel measure on $\Omega_C$, then there exists a $C$-full set $K$ such that
\begin{equation}\label{Normalized-form}
\frac{\lambda_\phi(K,\cdot)}{\lambda_\phi(K,\Omega_C)}=\frac{\mu}{\mu(\Omega_C)}\ \ \mbox{on}\ \ \Omega_C.
\end{equation}
Moreover, there are infinitely many such $K$. Formula \eqref{Normalized-form} should be understood in the sense of the limit, please see the proof in Theorem \ref{The-Main-Theorem-local-finite}.
\end{theorem}

\vskip 1mm

This paper is organized as follows. In Section \ref{Preliminaries}, we collect some preliminary materials on pseudo-cones. In Section \ref{Orlicz-integral Gauss curvature}, we introduce the Orlicz-Gauss image measure of pseudo-cones and its related properties. Section \ref{A Variational formula} is devoted to a variational formula for pseudo-cones under nonhomogeneous perturbations. In Section \ref{Existential results on compact sets}, we establish the uniform estimates and existence results on compact domains to the Orlicz-Gauss image problem for pseudo-cones. Section \ref{Proof of the main theorem} presents the proof of Theorem \ref{The-Main-Theorem}, and Section \ref{Proof of the main theorem local finite} that of Theorem \ref{The-Main-Theorem-local-finite}. Finally, we present the connection between the Orlicz-Gauss image problem for pseudo-cones and the spherical optimal transport in Section \ref{Section Spherical optimal transport}.

\section{Preliminaries}\label{Preliminaries}

This section collects some necessary background on pseudo-cones. We refer to \cite{Schneider-Pseudo_cones,Schneider-The_Gauss_image_problem_for_pseudo_cones} and the references therein for further details.

\subsection{$C$-determined sets and $C$-defined sets}

For $K\in ps(C)$, the copolar set of $K$ is defined by
\begin{equation*}
K^*=\{x\in\mathbb{R}^n:\langle x,y\rangle\leqslant-1, \forall\,y\in K\}\in ps(C^\circ).
\end{equation*}
Moreover, there holds $K^{**}=K$. The support function of $K$ is defined by
\begin{equation*}
h_K(v)=\sup_{x\in K}\langle v,x\rangle,\ v\in\operatorname{cl}\Omega_{C^\circ}.
\end{equation*}
Since $h_K\leqslant0$, we write $\overline{h}_K=-h_K$. Recall that the radial function of $K$ is defined by
\begin{equation*}
\rho_K(u)=\min\{r>0:ru\in K\},\ u\in \Omega_{C,K}.
\end{equation*}
The following relationships hold between the support function and the radial function of $K$:
\begin{align}
\rho_K(u)&=\frac{1}{\overline{h}_{K^*}(u)},\ \forall\,u\in\Omega_{C,K},\label{dual-relationship-support-radial}\\
\overline{h}_K(v)&=\inf_{u\in\Omega_C}|\langle u,v\rangle|\rho_K(u),\ \forall\,v\in\operatorname{cl}\Omega_{C^\circ},\label{Formula-need-OT}\\
\frac{1}{\rho_K(u)}&=\min_{v\in{\rm cl}\,\Omega_{C^\circ}}\frac{|\langle u,v\rangle|}{\overline{h}_K(v)},\ \forall\,u\in\Omega_C.\nonumber
\end{align}

Let $\omega\subset\Omega_{C^\circ}$ be a nonempty compact subset, then $K\in ps(C)$ is called $C$-determined by $\omega$ if
\begin{equation*}
K=C\cap\bigcap_{v\in\omega}H^-_K(v),
\end{equation*}
where $H^-_K(v)$ is the support halfspace of $K$ with the outer normal $v$. We denote by $\mathcal{K}(C,\omega)$ the class of all pseudo-cones that are $C$-determined by $\omega$. Given a positive continuous function $f$ on $\omega$ (written as $f\in\mathcal{C}^+(\omega)$), the Wulff shape associated with $(C,\omega,f)$ is defined by
\begin{equation*}
[f]=C\cap\bigcap_{v\in\omega}\{x\in\mathbb{R}^n:\langle x,v\rangle\leqslant-f(v)\}.
\end{equation*}
It is easy to check that $[f]\in\mathcal{K}(C,\omega)$.

Let $\eta\subset\Omega_C$ be a nonempty compact subset and $g\in\mathcal{C}^+(\eta)$, the convex hull associated with $(C,\eta,g)$ is defined by
\begin{equation*}
\langle g\rangle=\bigcap\{K\in ps(C):g(u)u\in K\ \mbox{for any}\ u\in\eta\}.
\end{equation*}
Clearly, $\langle g\rangle\in ps(C)$ and
\begin{equation*}
\rho_{\langle g\rangle}(u)\leqslant g(u)\ \ \mbox{for any}\ u\in\eta.
\end{equation*}
Note that $[g]\in ps(C^\circ)$, there holds
\begin{equation*}
[g]^*=\langle1/g\rangle.
\end{equation*}
Moreover, $K\in ps(C)$ is called $C$-defined by $\eta$ if $K=\langle\rho_K|_\eta\rangle$, where $\rho_K|_\eta$ is the restriction of the radial function $\rho_K$ to $\eta$. We denote by $\mathcal{K}^*(C,\eta)$ the class of all pseudo-cones that are $C$-defined by $\eta$. Then, there holds (see \cite{Schneider-The_Gauss_image_problem_for_pseudo_cones})
\begin{equation*}
K\in\mathcal{K}^*(C,\eta)\Leftrightarrow K^*\in\mathcal{K}(C^\circ,\eta).
\end{equation*}

\subsection{The radial Gauss map of pseudo-cones}

For $K\in ps(C)$, the radial map of $K$ is given by $r_K(u)=\rho_K(u)u$ for any $u\in\Omega_{C,K}$.
The Gauss image of $\sigma\subset\partial K$ is defined by
\begin{equation*}
\bm{\nu}_K(\sigma)=\{v\in\operatorname{cl}\Omega_{C^\circ}:\mbox{there exists a}\ x\in\sigma\ \mbox{such that}\ \langle x,v\rangle=h_K(v)\}.
\end{equation*}
The reverse Gauss image of $\omega\subset\operatorname{cl}\Omega_{C^\circ}$ is defined by
\begin{equation*}
\bm{x}_K(\omega)=\{x\in \partial K:\mbox{there exists a}\ v\in\omega\ \mbox{such that}\ \langle x,v\rangle=h_K(v)\}.
\end{equation*}
For $\eta\subset\Omega_{C,K}$, the radial Gauss image of $\eta$ is defined by $\bm{\alpha}_K(\eta)=\bm{\nu}_K(r_K(\eta))$. The reverse radial Gauss image of $\omega\subset{\rm cl}\,\Omega_{C^\circ}$ is defined by $\bm{\alpha}^*_K(\omega)=r^{-1}_K(\bm{x}_K(\omega))\subset\Omega_{C,K}$.

Let $\sigma_K\subset\partial K$ be the set of boundary points of $K$ where the unit outer normal is not unique, $\omega_K\subset\operatorname{cl}\Omega_{C^\circ}$ the set of $v\in\operatorname{cl}\Omega_{C^\circ}$ that serve as the unit outer normal at more than one point of $\partial K$, and $\eta_K$ the set of $u\in\Omega_{C,K}$ for which the unit outer normal at $r_K(u)$ is not unique. It is well-known from \cite{Schneider-book} that $\sigma_K$, $\omega_K$ and $\eta_K$ have $(n-1)$-dimensional Hausdorff measure zero, that is, $\mathcal{H}^{n-1}(\sigma_K)=\mathcal{H}^{n-1}(\omega_K)=\mathcal{H}^{n-1}(\eta_K)=0$.

For $x\in\partial K\setminus\sigma_K$, let $\nu_K(x)$ be the unique unit outer normal of $K$ at $x$. For $v\in{\rm cl}\,\Omega_{C^\circ}\setminus\omega_K$, let $x_K(v)$ be the unique boundary point of $K$ whose unit outer normal is $v$. Then, we define the radial Gauss map $\alpha_K:\Omega_{C,K}\setminus\eta_K\rightarrow\operatorname{cl}\Omega_{C^\circ}$ by
\begin{equation*}
\alpha_K(u)=\nu_K(r_K(u)),\ u\in\Omega_{C,K}\setminus\eta_K.
\end{equation*}
The reverse radial Gauss map $\alpha^*_K :\Omega_{C^\circ}\setminus\omega_K\rightarrow\Omega_{C,K}$ is defined by
\begin{equation*}
\alpha^*_K(v)=r^{-1}_K(x_K(v)),\ v\in\Omega_{C^\circ}\setminus\omega_K.
\end{equation*}
According to \cite{Schneider-The_Gauss_image_problem_for_pseudo_cones}, for $K\in ps(C)$, there holds
\begin{align}
\alpha^*_K(v)&=\alpha_{K^*}(v)\,\ \mbox{for}\,\ v\in\Omega_{C^\circ}\setminus\omega_K,\nonumber\\
\bm{\alpha}_K(\eta)\setminus\omega_K&=(\alpha^*_K)^{-1}(\eta)\,\ \mbox{for}\,\ \eta\subset\Omega_C.\label{Formula-alpha-image-reverse-map}
\end{align}

\subsection{Schneider selection theorem}

Let $K_i\in ps(C)$, $i\in\mathbb{N}$. We say that $\{K_i\}_{i=1}^{+\infty}$ converges to $K_0$ and write $K_i\rightarrow K_0$ as $i\rightarrow+\infty$, if there exists $t_0>0$ such that $K_i\cap t_0B^n\neq\emptyset$ for all $i\in\mathbb{N}$ and
\begin{equation*}
\lim_{i\to\infty}(K_i\cap tB^n)=K_0\cap tB^n\,\ \mbox{for all}\,\ t\geqslant t_0,
\end{equation*}
with respect to the Hausdorff metric. Here, $B^n=\{x\in\mathbb{R}^n:|x|\leqslant1\}$ is the unit ball. There holds the following equivalence \cite{Schneider-The_Gauss_image_problem_for_pseudo_cones}:
\begin{equation}\label{equivalence of convergence}
K_i\rightarrow K_0\ \Leftrightarrow K^*_i\rightarrow K^*_0.
\end{equation}

For $K\in ps(C)$, the distance of $K$ from the origin is denoted by $b(K)=\operatorname{dist}(o,K)$. Then, we have an important relationship as follows
\begin{equation}\label{relationship-support-function-distance-radial-function}
\overline{h}_K\leqslant b(K)\leqslant \rho_K.
\end{equation}
\begin{theorem}[Schneider selection theorem, \cite{Schneider-Pseudo_cones}]\label{Schneider selection theorem}
Let $K_i\in ps(C)$, $i\in\mathbb{N}$. If there exist two positive constants $a$ and $A$ such that
\begin{equation*}
a\leqslant b(K_i)\leqslant A,\ i\in\mathbb{N},
\end{equation*}
then there exists a subsequence of $\{K_i\}_{i=1}^{+\infty}$ converging to a $K\in ps(C)$.
\end{theorem}

\begin{lemma}[see \cite{Schneider-The_Gauss_image_problem_for_pseudo_cones}]\label{Radial-function-uniformly-convergence}
Let $K_i\in ps(C)$ converge to $K\in ps(C)$, and $\eta\subset\Omega_C$ be compact. Then
\begin{equation*}
\rho_{K_i}\rightarrow\rho_K\ \mbox{uniformly on}\ \eta, \ \mbox{as}\ i\to+\infty.
\end{equation*}
\end{lemma}

\begin{corollary}\label{Support-function-uniformly-convergence}
Let $K_i\in ps(C)$ converge to $K\in ps(C)$, and $\omega\subset\Omega_{C^\circ}$ be compact. Then
\begin{equation*}
\overline{h}_{K_i}\rightarrow\overline{h}_K\ \mbox{uniformly on}\ \omega, \ \mbox{as}\ i\to+\infty.
\end{equation*}
\end{corollary}
\begin{proof}
By \eqref{equivalence of convergence} and Lemma \ref{Radial-function-uniformly-convergence}, $\rho_{K_i^*}\rightarrow\rho_{K^*}$ uniformly on $\omega$, as $i\to+\infty$. Since $\omega$ is compact and $\rho_{K^*}$ is positive continuous, there is a positive constant $M$, such that
\begin{equation*}
\rho_{K^*},\rho_{K_i^*}>M\,\ \mbox{on}\,\ \omega,\,\ \mbox{for all}\,\ i\geqslant1.
\end{equation*}
And for any $\varepsilon>0$, there exists an $N>0$, such that
\begin{equation*}
|\rho_{K_i^*}-\rho_{K^*}|<M^2\varepsilon\,\ \mbox{on}\,\ \omega,\ \mbox{for all}\,\ i\geqslant N.
\end{equation*}
By \eqref{dual-relationship-support-radial}, we have
\begin{equation*}
\left|\frac{1}{\overline{h}_{K_i}}-\frac{1}{\overline{h}_{K}}\right|<M^2\varepsilon,
\end{equation*}
which reduces $\left|\bar{h}_{K_i}-\bar{h}_{K}\right|<M^2\varepsilon\,\bar{h}_{K_i}\bar{h}_{K}<\varepsilon$ on $\omega$.
\end{proof}

\section{Orlicz-Gauss image measure}\label{Orlicz-integral Gauss curvature}

This section is devoted to the introduction of the Orlicz-Gauss image measure of pseudo-cones and an investigation of its fundamental properties. Recall that $\lambda$ is a nonzero finite measure on spherical Lebesgue measurable subsets of $\operatorname{cl}\Omega_{C^\circ}$, which is zero on sets of Hausdorff dimension $n-2$, and that $\phi:(0,+\infty)\rightarrow(0,+\infty)$ is a continuous function.

\begin{definition}
Let $K\in ps(C)$, then the Orlicz-Gauss image measure of $K$ is defined by
\begin{equation*}
\lambda_\phi(K,\eta)=\int_{\bm{\alpha}_K(\eta)}\phi(\rho_K(\alpha^*_K(v)))\,\mathrm{d}\lambda(v),
\end{equation*}
for any Borel set $\eta\subset\Omega_C$. In particular, if $\phi\equiv1$, then $\lambda_\phi(K,\cdot)=\lambda(K,\cdot)$ is the Gauss image measure of $K$.
\end{definition}

We see that $\lambda_\phi(K,\cdot)$ is absolutely continuous with respect to $\lambda(K,\cdot)$, and there holds
\begin{equation*}
\mathrm{d}\lambda_\phi(K,u)=\phi(\rho_K(u))\mathrm{d}\lambda(K,u),
\end{equation*}
for $\lambda(K,\cdot)$-a.e. $u\in\Omega_C$. Moreover, since $\lambda_\phi(K,\cdot)$ is the pull-back of $\phi(\rho_K)\lambda$ via the radial Gauss map $\alpha_K$ (i.e., $\lambda_\phi(K,\cdot)$ is the push-forward of $\phi(\rho_K)\lambda$ via $\alpha^*_K$), we have the following important results:
\begin{proposition}
If $K\in ps(C)$ and $m$ is a Borel measure on $\Omega_C$, then the following terms are equivalent:\\
(i) $m=\lambda_\phi(K,\cdot)$;\\
(ii) For any bounded measurable function $f$ on $\Omega_C$, there holds
\begin{equation*}
\int_{\Omega_C}f(u)\,\mathrm{d}m(u)=\int_{\Omega_{C^\circ}}f(\alpha^*_K(v))\phi(\rho_K(\alpha^*_K(v)))\,\mathrm{d}\lambda(v);
\end{equation*}
(iii) For any bounded continuous function $f$ on $\Omega_C$, there holds
\begin{equation*}
\int_{\Omega_C}f(u)\,\mathrm{d}m(u)=\int_{\Omega_{C^\circ}}f(\alpha^*_K(v))\phi(\rho_K(\alpha^*_K(v)))\,\mathrm{d}\lambda(v).
\end{equation*}
\end{proposition}
\begin{proof}
{\bf Step 1: $(i)\Rightarrow(ii)$}. For any Borel set $\eta\subset\Omega_C$, by \eqref{Formula-alpha-image-reverse-map} we have
\[
v\in\bm{\alpha}_K(\eta)\setminus\omega_K\Leftrightarrow\alpha^*_K(v)\in\eta.
\]
It follows from \cite[Theorem 2.2.5]{Schneider-book} that $\omega_K$ can be covered by countably many sets of finite $(n-2)$-dimensional Hausdorff measure. This shows that $\lambda(\omega_K)=0$. Hence, there holds
\[
\chi_{\bm{\alpha}_K(\eta)}(v)=\chi_\eta(\alpha^*_K(v)) \ \ \text{for}\ \lambda\text{-a.e.}\ v\in\Omega_{C^\circ},
\]
where $\chi$ represents the characteristic function. Therefore, we have
\begin{align*}
\int_{\Omega_C}\chi_\eta(u)\,\mathrm{d}m(u)&=m(\eta)=\lambda_\phi(K,\eta)=\int_{\bm{\alpha}_K(\eta)}
\phi(\rho_K(\alpha^*_K(v)))\,\mathrm{d}\lambda(v)\\
&=\int_{\Omega_{C^\circ}}\chi_{\bm{\alpha}_K(\eta)}(v)\,\phi(\rho_K(\alpha^*_K(v)))\,\mathrm{d}\lambda(v)\\
&=\int_{\Omega_{C^\circ}}\chi_\eta(\alpha^*_K(v))\,\phi(\rho_K(\alpha^*_K(v)))\,\mathrm{d}\lambda(v).
\end{align*}
Then, for any simple function on $\Omega_C$
\begin{equation*}
\xi=\sum_{i=1}^ka_i\chi_{\eta_i}, \ \eta_i\subset\Omega_C,\ a_i\in\mathbb{R},
\end{equation*}
there holds
\begin{equation*}
\int_{\Omega_C}\xi(u)\,\mathrm{d}m(u)=\int_{\Omega_{C^\circ}}\xi(\alpha^*_K(v))\,\phi(\rho_K(\alpha^*_K(v)))\,\mathrm{d}\lambda(v).
\end{equation*}
Finally, since $f$ is a bounded measurable function, there exists a sequence of simple functions $\{\xi_i\}_{i=1}^{+\infty}$ converging uniformly to $f$ on $\Omega_C$. This concludes the desired (ii). Moreover, it is clear that (ii) implies (i) by taking the characteristic function of Borel sets.

{\bf Step 2: $(ii)\Rightarrow(iii)$}. It is also clear that (ii) implies (iii).

{\bf Step 3: $(iii)\Rightarrow(i)$}. From $(i)\Leftrightarrow(ii)$, we already know that for any bounded continuous function $f$ on $\Omega_C$, there holds
\[
\int_{\Omega_C}f(u)\,\mathrm{d}\lambda_\phi(K,u)=\int_{\Omega_{C^\circ}}f(\alpha^*_K(v))\phi(\rho_K(\alpha^*_K(v)))\,\mathrm{d}\lambda(v).
\]
Hence, we have
\[
\int_{\Omega_C}f(u)\,\mathrm{d}m(u)=\int_{\Omega_{C^\circ}}f(\alpha^*_K(v))\phi(\rho_K(\alpha^*_K(v)))\,\mathrm{d}\lambda(v)
=\int_{\Omega_C}f(u)\,\mathrm{d}\lambda_\phi(K,u).
\]
By the Riesz representation theorem, we conclude that $m=\lambda_\phi(K,\cdot)$.
\end{proof}

The following result now follows directly.
\begin{corollary}\label{Integral-curvature-pull-back-formula}
For any bounded measurable function $f$ on $\Omega_C$, there holds
\begin{equation*}
\int_{\Omega_C}f(u)\,\mathrm{d}\lambda_\phi(K,u)=\int_{\Omega_{C^\circ}}f(\alpha^*_K(v))\phi(\rho_K(\alpha^*_K(v)))\,\mathrm{d}\lambda(v).
\end{equation*}
\end{corollary}

Moreover, we have the following basic properties.
\begin{lemma}\label{Integral-curvature-concentrated-on-eta}
If $\eta\subset\Omega_{C}$ is compact and $K\in\mathcal{K}^*(C,\eta)$, then $\lambda_\phi(K,\cdot)$ is concentrated on $\eta$.
\end{lemma}
\begin{proof}
For $u\in\eta'=\Omega_C\backslash\eta$, since the unit outer normal of $K$ at $r_K(u)$ belongs to $\partial\Omega_{C^\circ}$, we have $\bm{\alpha}_{K}(\eta')\subset\partial\Omega_{C^\circ}$. This shows $\lambda(\bm{\alpha}_{K}(\eta'))=0$, thus
\begin{equation*}
\lambda_\phi(K,\eta')=\int_{\bm{\alpha}_K(\eta')}\phi(\rho_K(\alpha^*_K(v)))\,\mathrm{d}\lambda(v)=0.
\end{equation*}
This completes the proof of this lemma.
\end{proof}

\begin{lemma}\label{Orlicz-integral-of-full-set-is-finite}
If $K\in ps(C)$ is a $C$-full set, then $\lambda_\phi(K,\Omega_C)$ is finite.
\end{lemma}
\begin{proof}
Clearly, since $K$ is a $C$-full set and $\phi$ is continuous, so $\phi\circ\rho_K$ is bounded. Additionally, $\lambda$ is finite measure. Thus,
\begin{equation*}
\lambda_\phi(K,\Omega_C)=\int_{\bm{\alpha}_K(\Omega_C)}\phi(\rho_K(\alpha^*_K(v)))\,\mathrm{d}\lambda(v)
\end{equation*}
is finite.
\end{proof}

A sequence $\{K_i\}\subset ps(C)$ is called monotonic if either $K_i\subset K_{i+1}$ for all $i\in\mathbb{N}$, or $K_i\supset K_{i+1}$ for all $i\in\mathbb{N}$. For general Orlicz-Gauss image measures, we obtain the following important result concerning weak continuity.
\begin{proposition}\label{Weakly-convergence}
If a sequence of $C$-full sets $\{K_i\}_{i=1}^{+\infty}$ converges monotonically to a $C$-full set $K_0$, then $\lambda_\phi(K_i,\cdot)$ converges weakly to $\lambda_\phi(K_0,\cdot)$ on $\Omega_C$, as $i\to+\infty$.
\end{proposition}
\begin{proof}
Clearly, the radial function of a $C$-full set is continuous on $\operatorname{cl}\Omega_C$. Since $K_i\rightarrow K_0$ as $i\rightarrow+\infty$, by Lemma \ref{Radial-function-uniformly-convergence} we have $\rho_{K_i}(u)\rightarrow\rho_{K_0}(u)$ for all $u\in\Omega_C$. Since $K_i$ converges monotonically to $K_0$, there exists a $R>0$, such that convex bodies $K_i\cap RB^n$ monotonically converge to $K_0\cap RB^n$. This implies that $\rho_{K_i}(\bar{u})$ converges monotonically to $\rho_{K_0}(\bar{u})$ for any $\bar{u}\in\partial\Omega_C$. By the Dini theorem, $\rho_{K_i}$ converges uniformly to $\rho_{K_0}$ on $\operatorname{cl}\Omega_C$.

Let $f:\operatorname{cl}\Omega_C\rightarrow\mathbb{R}$ be a continuous function. From \cite[Lemma 4]{Schneider-The_Gauss_image_problem_for_pseudo_cones}, we have $\alpha^*_{K_i}(v)\rightarrow\alpha^*_{K_0}(v)$ for $\mathcal{H}^{n-1}$-a.e. $v\in\Omega_{C^\circ}$. Thus, for $\lambda$-a.e. $v\in\Omega_{C^\circ}$, there holds
\begin{equation*}
\lim_{i\rightarrow\infty}f(\alpha^*_{K_i}(v))\phi(\rho_{K_i}(\alpha^*_{K_i}(v)))=f(\alpha^*_{K_0}(v))\phi(\rho_{K_0}(\alpha^*_{K_0}(v))).
\end{equation*}
Due to the continuity of $f$, there exists a constant $M>0$, such that
\begin{equation*}
|f(\alpha^*_{K_i}(v))|\leqslant M,\ \mbox{for all}\ i\in\mathbb{N}.
\end{equation*}

Since $\rho_{K_i}$ converges uniformly to $\rho_{K_0}$ on $\operatorname{cl}\Omega_C$, we have $\rho_{K_i}<\rho_{K_0}+1$ for sufficiently large $i$. Thus, we conclude that
\begin{equation*}
\begin{aligned}
\int_{\Omega_{C^\circ}}|f(\alpha^*_{K_i}(v))\phi(\rho_{K_i}(\alpha^*_{K_i}(v)))|\,\mathrm{d}\lambda(v)&\leqslant M\int_{\Omega_{C^\circ}}\phi(\rho_{K_i}(\alpha^*_{K_i}(v)))\,\mathrm{d}\lambda(v)\\
&\leqslant M\int_{\Omega_{C^\circ}}\phi(\rho_{K_0}(\alpha^*_{K_i}(v))+1)\,\mathrm{d}\lambda(v)\\
&<C(\phi,K_0,M)<+\infty,
\end{aligned}
\end{equation*}
where $C$ is a constant depending on $\phi,K_0,M$. By the dominated convergence theorem and Corollary \ref{Integral-curvature-pull-back-formula}, we have
\begin{equation*}
\begin{aligned}
\lim_{i\rightarrow\infty}\int_{\Omega_C}f(u)\,\mathrm{d}\lambda_\phi(K_i,u)&=\lim_{i\rightarrow\infty}\int_{\Omega_{C^\circ}}f(\alpha^*_{K_i}(v))
\phi(\rho_{K_i}(\alpha^*_{K_i}(v)))\,\mathrm{d}\lambda(v)\\
&=\int_{\Omega_{C^\circ}}f(\alpha^*_{K_0}(v))\phi(\rho_{K_0}(\alpha^*_{K_0}(v)))\,\mathrm{d}\lambda(v)\\
&=\int_{\Omega_C}f(u)\,\mathrm{d}\lambda_\phi(K_0,u).
\end{aligned}
\end{equation*}
This completes the proof of this proposition.
\end{proof}

\section{A Variational formula}\label{A Variational formula}

This section establishes a variational formula that will play an essential role in the proof of our main result.
Denote by $\mathcal{J}_{\delta,a,b}$ the set of continuously differentiable functions $\varphi\colon(\delta,+\infty)\rightarrow(a,b)$ satisfying
\begin{equation*}
\begin{cases}
\varphi'>0,\\
\displaystyle\lim_{t\to\delta^+}\varphi(t)=a,\\
\displaystyle\lim_{t\to+\infty}\varphi(t)=b,
\end{cases}
\end{equation*}
where $a\in\mathbb{R}$, $b\in\mathbb{R}\cup\{+\infty\}$, $\delta\geqslant0$, and $a<b$. Let $\eta\subset\Omega_C$ be compact, $f_0$ be a continuous function on $\eta$ satisfying $f_0>\delta$ (written as $f_0\in\mathcal{C}_\delta(\eta)$), and $g$ be a continuous function on $\eta$ (written as $g\in\mathcal{C}(\eta)$). For $\varphi\in\mathcal{J}_{\delta,a,b}$, we define a family of $\varphi$-perturbations as follows:
\begin{equation}\label{varphi-perturbations}
\varphi(f_t(u))=\varphi(f_0(u))+tg(u)\,\ \mbox{for any}\,\ u\in\eta,
\end{equation}
where $f_t\in\mathcal{C}_\delta(\eta)$ for sufficiently small $t$ (i.e., $t\to0$).

Using the $\varphi$-perturbations, the variational formulas for convex bodies are first studied by Gardner et al. in \cite{Gardner-General_volumes_in_the_Orlicz_Brunn_Minkowski_theory}. Here, we give the corresponding variational formulas for pseudo-cones. Notably, our method is significantly more convenient.

\begin{lemma}\label{The-variational-lemma}
Let $\eta\subset\Omega_C$ be compact, $\varphi\in\mathcal{J}_{\delta,a,b}$, $f_0\in\mathcal{C}_\delta(\eta)$ and $g\in\mathcal{C}(\eta)$. For the $\varphi$-perturbations $f_t$ defined by \eqref{varphi-perturbations} and any $v\in\Omega_{C^\circ}\setminus\omega_{[f_0]}$, there holds
\begin{equation*}
\lim_{t\rightarrow 0}\frac{\log\overline{h}_{\langle f_t\rangle}(v)-\log\overline{h}_{\langle f_0\rangle}(v)}{t}
=\frac{g(\alpha^*_{\langle f_0\rangle}(v))}{f_0(\alpha^*_{\langle f_0\rangle}(v))\varphi'(f_0(\alpha^*_{\langle f_0\rangle}(v)))}.
\end{equation*}
Moreover, there exist $\varepsilon>0$ and $M>0$ (depending on $g,f_0,\varphi$), such that
\begin{equation*}
|\log\overline{h}_{\langle f_t\rangle}(v)-\log\overline{h}_{\langle f_0\rangle}(v)|\leqslant M|t|\,\ \mbox{for all}\,\ t\in(-\varepsilon,\varepsilon)\ \mbox{and}\ v\in\Omega_{C^\circ}.
\end{equation*}
\end{lemma}
\begin{proof}
Let $H_u(t):=\log f_t(u)=\log\varphi^{-1}(\varphi(f_0(u))+tg(u))$, $u\in\eta$, then
\begin{equation*}
H'_u(t)=\frac{g(u)}{f_t(u)\varphi'(f_t(u))}.
\end{equation*}
Using the Taylor's formula, we have
\begin{equation*}
H_u(t)=H_u(0)+tH'_u(0)+o(u,t),
\end{equation*}
where $o(u,t)$ is an infinitesimal of higher order.

Next, we prove that $H'_u(t)$ converges uniformly to $H'_u(0)$ on $\eta$, as $t\to0$. Note that $tg(u)$ converges uniformly to $0$ on $\eta$, as $t\to0$, thus $\varphi(f_0(u))+tg(u)\to\varphi(f_0(u))$ uniformly on $\eta$, as $t\to0$. Since $\varphi(f_0(u))$ is continuous on $\eta$, $\varphi(f_0(u))+tg(u)$ is uniformly bounded on $\eta$. Combining with the continuity of $\varphi^{-1}$, we conclude that
\[
\varphi^{-1}(\varphi(f_0(u))+tg(u))\to\varphi^{-1}(\varphi(f_0(u)))
\]
uniformly on $\eta$, as $t\to0$. That is, $f_t(u)\to f_0(u)$ uniformly on $\eta$, as $t\to0$. Again, due to the continuity of $\varphi'$, we have $\varphi'(f_t(u))\to\varphi'(f_0(u))$ uniformly on $\eta$, as $t\to0$. Let $M_1=\max_{u\in\eta}|g(u)|$. Choose a sufficiently small $\delta_0>0$ depending on $M_1$, such that for all $t\in(-\delta_0,\delta_0)$ and $u\in\eta$, there holds
\begin{equation*}
f_t(u)\geqslant\varphi^{-1}(\varphi(f_0(u))-\delta_0M_1)=:M_2(u)>0.
\end{equation*}
Note that $M_2(u)$ is continuous, $f_t(u)$ and $\varphi'(f_t(u))$ are uniformly bounded from below on $\eta$. This gives
\begin{equation*}
\frac{g(u)}{f_t(u)\varphi'(f_t(u))}\to\frac{g(u)}{f_0(u)\varphi'(f_0(u))}\ \mbox{uniformly on}\ \eta,\ \mbox{as}\ t\to0.
\end{equation*}

Combining the Lagrange's mean value theorem, we have
\begin{equation*}
\begin{aligned}
\left|\frac{o(u,t)}{t}\right|&=\left|\frac{H_u(t)-H_u(0)}{t}-H'_u(0)\right|\\
&=|H'_u(\theta t)-H'_u(0)|\to0\,\ \mbox{uniformly on}\,\ \eta,\ \mbox{as}\ t\to0,
\end{aligned}
\end{equation*}
where $\theta\in(0,1)$. Using the variational formula with higher-order infinitesimals in
\cite[Lemma 11]{Schneider-The_Gauss_image_problem_for_pseudo_cones} (it can also be derived from
\cite[Lemma 9]{Schneider-A_weighted_Minkowski_theorem_for_pseudo_cones} or
\cite[Lemma 5.4]{Li-Zhu-Ye-The_dual_Minkowski_problem_for_unbounded_closed_convex_sets}), we have
\begin{align*}
\lim_{t\rightarrow 0}\frac{\log\overline{h}_{\langle f_t\rangle}(v)-\log\overline{h}_{\langle f_0\rangle}(v)}{t}&=
H'_{\alpha^*_{\langle f_0\rangle}(v)}(0)\\
&=\frac{g(\alpha^*_{\langle f_0\rangle}(v))}{f_0(\alpha^*_{\langle f_0\rangle}(v))\varphi'(f_0(\alpha^*_{\langle f_0\rangle}(v)))},
\end{align*}
and there exist $\varepsilon>0$ and $M>0$ (depending on $g,f_0,\varphi$), such that
\[
|\log\overline{h}_{\langle f_t\rangle}(v)-\log\overline{h}_{\langle f_0\rangle}(v)|\leqslant M|t|,
\]
for all $v\in\Omega_{C^\circ}$ and $t\in(-\varepsilon,\varepsilon)$.
\end{proof}

\begin{proposition}\label{The-variational-formula}
Let $\eta\subset\Omega_C$ be compact, $\varphi\in\mathcal{J}_{\delta,a,b}$, $K\in\mathcal{K}^*(C,\eta)$ and $g\in\mathcal{C}(\eta)$. If $b(K)>\delta$ and $\varphi(f_t(u))=\varphi(\rho_K(u))+tg(u)$ for all $u\in\eta$, then
\begin{equation*}
\frac{d}{dt}\bigg|_{t=0}\int_{\Omega_{C^\circ}}\log\overline{h}_{\langle f_t\rangle}(v)\,\mathrm{d}\lambda(v)=\int_\eta g(u)d\lambda_\psi(K,u),
\end{equation*}
where $\psi(t)=\frac{1}{t\varphi'(t)}$.
\end{proposition}
\begin{proof}
Let $f_0$ be the restriction of $\rho_K$ on $\eta$ (written as $f_0=\rho_K|\eta$), then $K=\langle f_0\rangle$. From the dominated convergence theorem, Corollary \ref{Integral-curvature-pull-back-formula}, Lemma \ref{Integral-curvature-concentrated-on-eta} and Lemma \ref{The-variational-lemma}, we have
\begin{equation*}
\begin{aligned}
\frac{d}{dt}\bigg|_{t=0}\int_{\Omega_{C^\circ}}\log\overline{h}_{\langle f_t\rangle}(v)\,\mathrm{d}\lambda(v)
&=\lim_{t\to0}\int_{\Omega_{C^\circ}}\frac{\log\overline{h}_{\langle f_t\rangle}(v)-\log\overline{h}_K(v)}{t}\,\mathrm{d}\lambda(v)\\
&=\int_{\Omega_{C^\circ}}\frac{g(\alpha^*_K(v))}{\rho_K(\alpha^*_K(v))\varphi'(\rho_K(\alpha^*_K(v)))}\,\mathrm{d}\lambda(v)\\
&=\int_{\Omega_{C^\circ}}g(\alpha^*_K(v))\psi(\rho_K(\alpha^*_K(v)))\,\mathrm{d}\lambda(v)\\
&=\int_{\Omega_{C}}g(u)\,\mathrm{d}\lambda_\psi(K,u)\\
&=\int_\eta g(u)\,\mathrm{d}\lambda_\psi(K,u).
\end{aligned}
\end{equation*}
This establishes the variational formula.
\end{proof}

\section{Existential results on compact sets}\label{Existential results on compact sets}

This section proves the existence of solutions to the Orlicz-Gauss image problem on compact sets. First, we define the entropy of $K\in ps(C)$ by
\begin{equation*}
\mathcal{E}(K)=-\int_{\Omega_{C^\circ}}\log\overline{h}_K(v)\,\mathrm{d}\lambda(v).
\end{equation*}
Then, we have the following basic yet important estimate.
\begin{lemma}[The lower bound estimate]\label{The-lower-bound-estimate}
For $K\in ps(C)$, if $\mathcal{E}(K)=1$, then there exists a constant $\delta>0$ such that $b(K)\geqslant\delta$.
\end{lemma}
\begin{proof}
Using \eqref{relationship-support-function-distance-radial-function}, we have
\begin{equation*}
-1=\int_{\Omega_{C^\circ}}\log\overline{h}_K(v)\,\mathrm{d}\lambda(v)\leqslant\lambda(\Omega_{C^\circ})\log b(K),
\end{equation*}
then
\begin{equation*}
b(K)\geqslant\exp\left(-\frac{1}{\lambda(\Omega_{C^\circ})}\right)=:\delta.
\end{equation*}
This yields the desired estimate.
\end{proof}

The following result shows that the entropies of $C$-defined sets are continuous.
\begin{lemma}\label{Continuity-entropy}
Let $\eta\subset\Omega_C$ be compact and $\{K_i\}_{i=0}^{+\infty}\subset\mathcal{K}^*(C,\eta)$. If $K_i$ converges to $K_0$ as $i\to+\infty$, then
\begin{equation*}
\lim_{i\to+\infty}\mathcal{E}(K_i)=\mathcal{E}(K_0).
\end{equation*}
\end{lemma}
\begin{proof}
Since $K_i\to K_0$, there exists a $R>0$ such that $K_i\cap RB^n\to K_0\cap RB^n$ in the Hausdorff metric. This implies that $\operatorname{dist}(o,K_i\cap RB^n)\to\operatorname{dist}(o,K_0\cap RB^n)$, i.e., $b(K_i)\to b(K_0)$. Thus, there exists $N>0$ such that $b(K_i)\leqslant b(K_0)+1$ for any $i>N$. By Corollary \ref{Support-function-uniformly-convergence}, we know that $\overline{h}_{K_i}(v)\to\overline{h}_{K_0}(v)$ for any $v\in\Omega_{C^\circ}$. Then, for any $v\in\Omega_{C^\circ}$ and $i>N$, we have
\begin{equation*}
\overline{h}_{K_i}(v)\leqslant b(K_i)\leqslant b(K_0)+1.
\end{equation*}
Let $A=\max\{b(K_1),\cdots,b(K_N),b(K_0)+1\}$, then $\overline{h}_{K_i}(v)\leqslant A$ for any $v\in\Omega_{C^\circ}$ and $i\in\mathbb{N}$.

On the other hand, since $\eta\subset\Omega_C$ is compact and $\{K_i\}_{i=0}^{+\infty}\subset\mathcal{K}^*(C,\eta)$, there exists a positive constant $a$ such that $\overline{h}_{K_i}\geqslant a$. Therefore, we have
\begin{equation*}
|\log\overline{h}_{K_i}(v)|\leqslant\max\{|\log a|,|\log A|\}\,\ \mbox{for any}\,\ v\in\Omega_{C^\circ}\,\ \mbox{and}\,\ i\in\mathbb{N}.
\end{equation*}
Finally, the dominated convergence theorem leads to the desired conclusion.
\end{proof}

Now, by the variational method, we obtain existence results on compact domains.
\begin{theorem}\label{Theorem on compact set}
Let $\eta\subset\Omega_C$ be compact and $\mu$ be a nonzero finite Borel measure on $\eta$. Then, there exists a $K\in\mathcal{K}^*(C,\eta)$ such that
\begin{equation*}
\frac{\lambda_{\phi}(K,\cdot)}{\lambda_{\phi}(K,\eta)}=\frac{\mu}{\mu(\eta)}.
\end{equation*}
\end{theorem}
\begin{proof}
By Lemma \ref{The-lower-bound-estimate}, there exists a $\delta>0$, such that $b(K)>\delta$ for any $K\in ps(C)$ with $\mathcal{E}(K)=1$. We define the function
\begin{equation*}
\varphi(t)=\int_\delta^t\frac{1}{s\phi(s)}\,ds\ \ \mbox{for}\ \ t>\delta,
\end{equation*}
then $\varphi\in\mathcal{J}_{\delta,0,b}$, where $b\in(0,+\infty]$. We define the functional
\begin{equation*}
\mathcal{F}(f)=\frac{1}{\mu(\eta)}\int_{\eta}\varphi(f(u))\,\mathrm{d}\mu(u),\ f\in\mathcal{C}_\delta(\eta).
\end{equation*}
In particular, we write $\mathcal{F}(K)=\mathcal{F}(\rho_K|\eta)$ for any $K\in\mathcal{K}^*(C,\eta)$ with $\mathcal{E}(K)=1$. Clearly, $\mathcal{F}$ is continuous on $\mathcal{C}_\delta(\eta)$. That is, if $f_i\in\mathcal{C}_\delta(\eta)$ converges uniformly to $f\in\mathcal{C}_\delta(\eta)$, then $\mathcal{F}(f_i)$ converges to $\mathcal{F}(f)$.

Consider the infimum
\begin{equation*}
\gamma=\inf\{\mathcal{F}(f):f\in\mathcal{C}_\delta(\eta)\,\ \text{and}\,\ \mathcal{E}(\langle f\rangle)=1\}\geqslant0.
\end{equation*}
It can be found from $\rho_{\langle f\rangle}\leqslant f$ that this infimum can only be attained by the radial function of some pseudo-cones that are $C$-defined by $\eta$. This shows that
\begin{equation*}
\inf\{\mathcal{F}(K):K\in\mathcal{K}^*(C,\eta)\,\ \mbox{and}\,\ \mathcal{E}(K)=1\}=\gamma.
\end{equation*}

Let $\left\{K_i\right\}_{i=1}^{+\infty}\subset\mathcal{K}^*(C,\eta)$ be a minimizing sequence, that is, $\mathcal{E}(K_i)=1$ and
\begin{equation*}
\lim_{i\to\infty}\mathcal{F}(K_i)=\gamma.
\end{equation*}
Thus, there exists a $N>0$, such that $\mathcal{F}(K_i)<\gamma+1$ for all $i\geqslant N$. Since $b_{K_i}\leqslant\rho_{K_i}$, we have
\begin{equation*}
\varphi(b_{K_i})=\frac{1}{\mu(\eta)}\int_{\eta}\varphi(b_{K_i})\,\mathrm{d}\mu(u)\leqslant\mathcal{F}(K_i).
\end{equation*}
We choose $c=\max\{\mathcal{F}(K_1),\cdots,\mathcal{F}(K_N),\gamma+1\}$, then $\varphi(b_{K_i})\leqslant c$. Thus, we have
\begin{equation*}
\delta\leqslant b_{K_i}\leqslant\varphi^{-1}(c),\ i=1,2,\cdots.
\end{equation*}
By Theorem \ref{Schneider selection theorem}, the sequence $\{K_i\}_{i=1}^{+\infty}$ has a subsequence converging to a $K\in ps(C)$. Clearly, we also have $K\in\mathcal{K}^*(C,\eta)$. W.l.o.g., we assume that $\{K_i\}_{i=1}^{+\infty}$ converges to $K$. By Lemma \ref{Radial-function-uniformly-convergence}, $\rho_{K_i}$ converges uniformly to $\rho_K$ on $\eta$. Thus, $\mathcal{F}(K_i)$ converges to $\mathcal{F}(K)$. Therefore, $K\in\mathcal{K}^*(C,\eta)$ satisfies $\mathcal{E}(K)=1$ (by Lemma \ref{Continuity-entropy}) and
\begin{equation*}
\mathcal{F}(K)=\inf\{\mathcal{F}(f):f\in\mathcal{C}_\delta(\eta)\,\ \mbox{and}\,\ \mathcal{E}(\langle f\rangle)=1\}.
\end{equation*}

Given a $g\in\mathcal{C}(\eta)$, we define $f_t$ by $\varphi(f_t(u))=\varphi(\rho_K(u))+tg(u)$ for all $u\in\eta$. Consider the Lagrange function $\mathcal{L}(t,\tau)=\mathcal{F}(f_t)+\tau(\mathcal{E}(\langle f_t\rangle)-1)$, then
\begin{equation*}
\frac{\partial}{\partial t}\bigg|_{t=0}\mathcal{L}(t,\tau)=0.
\end{equation*}
By Proposition \ref{The-variational-formula}, we have
\begin{equation*}
\begin{aligned}
0&=\frac{\partial}{\partial t}\bigg|_{t=0}(\mathcal{F}(f_t)+\tau(\mathcal{E}(\langle f_t\rangle)-1))\\
&=\frac{\partial}{\partial t}\bigg|_{t=0}\left(\frac{1}{\mu(\eta)}\int_\eta\varphi(\rho_K(u))+tg(u)\,\mathrm{d}\mu(u)\right)
+\tau\frac{\partial}{\partial t}\bigg|_{t=0}\mathcal{E}(\langle f_t\rangle)\\
&=\frac{1}{\mu(\eta)}\int_\eta g(u)\,\mathrm{d}\mu(u)-\tau\int_\eta g(u)\,\mathrm{d}\lambda_\phi(K,u),
\end{aligned}
\end{equation*}
That is,
\begin{equation*}
  \frac{1}{\mu(\eta)}\int_\eta g(u)\,\mathrm{d}\mu(u)=\tau\int_\eta g(u)\,\mathrm{d}\lambda_\phi(K,u).
\end{equation*}
By the Riesz representation theorem, we conclude that $\mu=\mu(\eta)\tau \lambda_\phi(K,\cdot)$. Clearly, we can conclude that $\tau=1/\lambda_\phi(K,\eta)$.
\end{proof}

\section{Proof of Theorem \ref{The-Main-Theorem}}\label{Proof of the main theorem}

This section gives the proof of Theorem \ref{The-Main-Theorem}. To this end, we need the following key tool.
\begin{proposition}\label{Double-cone-Lemma}
Let $C$ and $\Gamma$ be two $n$-dimensional pointed closed convex cone. If $C\subset\Gamma$ and $K\in\mathcal{K}^*(\Gamma,\mathbb{S}^{n-1}\cap C)$, then $K\cap C$ is a $C$-full set. Moreover, we have
\begin{equation*}
\lambda_\phi(K,\cdot)=\lambda_\phi(K\cap C,\cdot)\,\ \mbox{on}\,\ \Omega_C.
\end{equation*}
\end{proposition}
\begin{proof}
Clearly, $K\cap C$ is a nonempty closed convex set and $C\setminus(K\cap C)$ is bounded, so $K\cap C$ is a $C$-full set. Note that $\rho_K=\rho_{K\cap C}$ on $\operatorname{cl}\Omega_C$, then we have
\begin{equation*}
\bm{\alpha}_K(\eta)=\bm{\alpha}_{K\cap C}(\eta)\subset\operatorname{cl}\Omega_{\Gamma^\circ}\,\ \mbox{for all}\,\ \eta\subset\Omega_C,
\end{equation*}
and
\begin{equation*}
\alpha^*_K(v)=\alpha^*_{K\cap C}(v)\in\operatorname{cl}\Omega_C\,\ \mbox{for any}\,\ v\in\operatorname{cl}\Omega_{\Gamma^\circ}\setminus\omega_K.
\end{equation*}
This yields that
\begin{equation*}
\begin{aligned}
\lambda_\phi(K,\eta)&=\int_{\bm{\alpha}_K(\eta)}\phi(\rho_K(\alpha^*_K(v)))\,\mathrm{d}\lambda(v)\\
&=\int_{\bm{\alpha}_{K\cap C}(\eta)}\phi(\rho_{K\cap C}(\alpha^*_{K\cap C}(v)))\,\mathrm{d}\lambda(v)\\
&=\lambda_\phi(K\cap C,\eta),
\end{aligned}
\end{equation*}
for all $\eta\subset\Omega_C$.
\end{proof}

For $K\in ps(C)$, it is easy to see that the Orlicz-Gauss image measure $\lambda_\phi(K,\cdot)$ can be defined on $\operatorname{cl}\Omega_C$. However, in the above Proposition \ref{Double-cone-Lemma}, the equality $\lambda_\phi(K,\cdot)=\lambda_\phi(K\cap C,\cdot)$ does not hold on $\operatorname{cl}\Omega_C$. Consider the following counterexample in the planar case.

\begin{figure}[http]
  \centering
  \includegraphics[width=8cm]{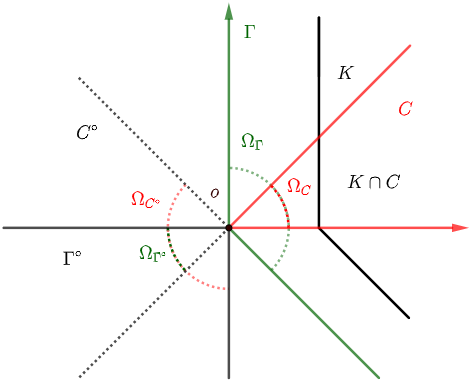}
  \caption{An counterexample: $K$ is $\Gamma$-defined by $\operatorname{cl}\Omega_C$.}\label{Gauss_image_measure_counterexample}
\end{figure}
\noindent As shown in Figure \ref{Gauss_image_measure_counterexample}, we choose $\lambda=\mathcal{H}^1$, $\phi=1$, and
\begin{align*}
\Gamma&=\{(x,y)\in\mathbb{R}^2:x\geqslant0,y\geqslant-x\},\\
C&=\{(x,y)\in\mathbb{R}^2:x\geqslant0,0\leqslant y\leqslant x\},\\
K&=\{(x,y)\in\mathbb{R}^2:x\geqslant2,y\geqslant2-x\}.
\end{align*}
Then, $\bm{\alpha}_K(\operatorname{cl}\Omega_C)=\operatorname{cl}\Omega_{\Gamma^\circ}$ and $\bm{\alpha}_{K\cap C}(\operatorname{cl}\Omega_C)=\operatorname{cl}\Omega_{C^\circ}$. Hence,
\begin{align*}
&\lambda_\phi(K,\operatorname{cl}\Omega_C)=\int_{\operatorname{cl}\Omega_{\Gamma^\circ}}\,\mathrm{d}\mathcal{H}^1
=\mathcal{H}^1(\operatorname{cl}\Omega_{\Gamma^\circ}),\\
&\lambda_\phi(K\cap C,\operatorname{cl}\Omega_C)=\int_{\operatorname{cl}\Omega_{C^\circ}}\,\mathrm{d}\mathcal{H}^1
=\mathcal{H}^1(\operatorname{cl}\Omega_{C^\circ}).
\end{align*}

Now, using the restrictive technique, we give the proof of Theorem \ref{The-Main-Theorem} as follows.
\begin{proof}[Proof of Theorem \ref{The-Main-Theorem}]
Let $\mu$ be a nonzero Borel measure on $\Omega_C$. If $\mu$ is finite, then we define the measure
\begin{equation*}
\widetilde{\mu}(\eta)=\mu(\eta\cap\Omega_C)\ \ \mbox{for any}\ \,\eta\subset{\rm cl}\,\Omega_C,
\end{equation*}
which is a nonzero finite Borel measure on ${\rm cl}\,\Omega_C$. For the cone $C$, we can choose a larger pointed closed convex cone $\Gamma\supset C$ (Since $C$ is pointed, we can certainly find such a $\Gamma$). Denote $\eta={\rm cl}\,\Omega_C$, then $\eta\subset\Omega_\Gamma$ is compact. Then, by Theorem \ref{Theorem on compact set}, there exists a $L\in\mathcal{K}^*(\Gamma,\eta)$, such that
\begin{equation*}
\frac{\lambda_\phi(L,\cdot)}{\lambda_\phi(L,\eta)}=\frac{\widetilde{\mu}}{\widetilde{\mu}(\eta)}\,\ \mbox{on}\,\ \eta.
\end{equation*}
Let $K=L\cap C$, then from Proposition \ref{Double-cone-Lemma} we know that $K$ is a $C$-full set and
\begin{equation*}
\lambda_\phi(L,\cdot)=\lambda_\phi(K,\cdot)\,\ \mbox{on}\,\ \Omega_C.
\end{equation*}
Therefore, we have
\begin{equation*}
\mu=\widetilde{\mu}=\frac{\widetilde{\mu}(\eta)}{\lambda_\phi(L,\eta)}\lambda_\phi(K,\cdot)=
\frac{\mu(\Omega_C)}{\lambda_\phi(L,\eta)}\lambda_\phi(K,\cdot)=:\beta\lambda_\phi(K,\cdot)\,\ \mbox{on}\,\ \Omega_C.
\end{equation*}
Conversely, if there is a $C$-full set $K$ and a constant $\beta>0$ such that $\beta\lambda_\phi(K,\cdot)=\mu$, then by Lemma \ref{Orlicz-integral-of-full-set-is-finite} we know that $\lambda_\phi(K,\cdot)$ is finite. Thus, $\mu$ is finite.

In addition, we note that the constant $\beta$ in the above proof depends on the choice of the larger cone $\Gamma$. Hence, generally speaking, such $\beta$ and $K$ are not unique. That is, if there exist $K_1,K_2\in ps(C)$ and two positive constants $\beta_1\neq\beta_2$ such that
\begin{equation*}
\beta_1\lambda_\phi(K_1,\cdot)=\beta_2\lambda_\phi(K_2,\cdot)\,\ \mbox{on}\,\ \Omega_C,
\end{equation*}
then $K_1\neq K_2$. Due to the freedom in choosing the greater cone, such constant $c$ admits infinitely many distinct values. Consequently, we derive infinitely many solutions to the Orlicz-Gauss image problem for pseudo-cones.
\end{proof}

The restrictive technique, which does not rely on uniform estimates, is different from Schneider's approximation approach. Consequently, we strengthen Theorem \ref{Schneider-theorem} by a constant factor. Specifically, we prove Corollary \ref{The-corollary-1} as follows.
\begin{proof}[Proof of Corollary \ref{The-corollary-1}]
Let $\mu$ and $\lambda$ are probability measures on $\Omega_C$ and $\operatorname{cl}\Omega_{C^\circ}$ respectively. For $\phi=1$, by Theorem \ref{The-Main-Theorem}, there exists a $C$-full set $K$ and $\beta>0$, such that
\[
\beta\lambda(K,\cdot)=\mu.
\]
Thus, we have
\[
\beta=\frac{\mu(\Omega_C)}{\lambda(K,\Omega_C)}=\frac{\mu(\Omega_C)}{\lambda(\bm{\alpha}_K(\Omega_C))}\geqslant
\frac{\mu(\Omega_C)}{\lambda(\operatorname{cl}\Omega_{C^\circ})}=1,
\]
where $\bm{\alpha}_K(\Omega_C)\subsetneqq\operatorname{cl}\Omega_{C^\circ}$ (since $K$ is a $C$-full set). That is, there exists a $C$-full set $K$ and a constant $\beta'\in(0,1]$ such that
\[
(\alpha^*_K)_\sharp\,\lambda=\beta'\mu.
\]
\end{proof}

For the $L_p$ case, we give the proof of Corollary \ref{The-corollary-2}.
\begin{proof}[Proof of Corollary \ref{The-corollary-2}]
Let $\mu$ be a nonzero Borel measure on $\Omega_C$ and $0\neq p\in\mathbb{R}$. We choose $\phi(t)=t^p$. Clearly, if there exists a $C$-full set $K$ such that $\lambda_p(K,\cdot)=\mu$, then $\mu$ is finite.

On the other hand, if $\mu$ is finite, by Theorem \ref{The-Main-Theorem}, there exists a $C$-full set $K$ and $\beta>0$, such that
\[
\beta\lambda_p(K,\cdot)=\mu.
\]
For any Borel set $\eta\subset\Omega_C$, we have
\begin{align*}
\lambda_p(tK,\eta)&=\int_{\bm{\alpha}_{tK}(\eta)}\rho_{tK}^p(\alpha^*_{tK}(v))\,\mathrm{d}\lambda(v)\\
&=\int_{\bm{\alpha}_K(\eta)}t^p\rho_K^p(\alpha^*_K(v))\,\mathrm{d}\lambda(v)\\
&=t^p\lambda_p(K,\eta).
\end{align*}
Therefore, $\widetilde{K}=\beta^\frac{1}{p}K$ is a $C$-full set and satisfies $\lambda_p(\widetilde{K},\cdot)=\mu$.
\end{proof}

Finally, as Schneider pointed out in \cite{Schneider-The_Gauss_image_problem_for_pseudo_cones}, the Gauss image problem for pseudo-cones is simplified by considering measures $\mu$ on $\Omega_C$ rather than on its closure. Based on the previous counterexample regarding Proposition \ref{Double-cone-Lemma}, we cannot deal with measures $\mu$ on $\operatorname{cl}\Omega_C$. Thus, the generalization of all the present results to measures on $\operatorname{cl}\Omega_C$ remains open.

\section{The case of infinite measure and the proof of Theorem \ref{The-Main-Theorem-local-finite}}
\label{Proof of the main theorem local finite}

A measure $\mu$ on $\Omega_C$ is called locally finite if $\mu(\eta)<+\infty$ for any compact $\eta\subset\Omega_C$. In this section, we give the proof of Theorem \ref{The-Main-Theorem-local-finite}.

\begin{proof}[Proof of Theorem \ref{The-Main-Theorem-local-finite}]
We divide two cases as follows:

{\bf Case 1.} If $\mu$ is finite, i.e., $\mu(\Omega_C)<+\infty$, then by Theorem \ref{The-Main-Theorem} there is a $C$-full set $K$ and a constant $\beta>0$ such that $\beta\lambda_\phi(K,\cdot)=\mu$. By Lemma \ref{Orlicz-integral-of-full-set-is-finite}, $\lambda_\phi(K,\Omega_C)<+\infty$. Thus, we have
\begin{equation*}
\beta=\frac{\mu(\Omega_C)}{\lambda_\phi(K,\Omega_C)}.
\end{equation*}
This derives the desired result.

{\bf Case 2.} If $\mu(\Omega_C)=+\infty$, since $\mu$ is a locally finite, $\mu(\eta)<+\infty$ for any compact $\eta\subset\Omega_C$. Then,
\begin{equation*}
\frac{\mu(\eta)}{\mu(\Omega_C)}=0.
\end{equation*}
Choose a point $v_0\in\operatorname{int}\Omega_{C^\circ}$, we set
\[
K=H^-(v_0,-1)\cap C,
\]
where $H^-(v_0,-1)=\{x\in\mathbb{R}^n:\langle x,v_0\rangle\leqslant-1\}$. Then, $K$ is a $C$-full set and $\bm{\alpha}_K(\Omega_C)=\{v_0\}$. This shows that $\lambda_\phi(K,\Omega_C)=0$.

We choose a sequence of spherical closed convex sets $\{\eta_i\}_{i=1}^{+\infty}$ in $\Omega_C$ such that
\begin{equation*}
\eta_i\subset\eta_{i+1}\ \ \mbox{and}\ \ \bigcup_{i=1}^{+\infty}\eta_i=\Omega_C.
\end{equation*}
For each $\eta_i$, we construct the set
\begin{equation*}
K_i=\bigcap\{L\in ps(C):r_K(\eta_i)\cup\{r_K(u)+u:u\in\partial\Omega_C\}\subset L\},
\end{equation*}
then $K_i$ is also a $C$-full set. Moreover, $K_i\subset K_{i+1}$ for all $i\in\mathbb{N}$ and $K_i\rightarrow K$ as $i\rightarrow\infty$.

For any compact $\eta\subset\Omega_C$, there is $N>0$ such that $\eta\subset\eta_i$ for any $i>N$. Then, for $i>N$, we have $\bm{\alpha}_{K_i}(\eta)=\{v_0\}$. Hence, $\lambda_\phi(K_i,\eta)=0$. But, it is easy to check that $\lambda_\phi(K_i,\Omega_C)>0$. Thus, we conclude that
\begin{equation*}
\frac{\lambda_\phi(K_i,\eta)}{\lambda_\phi(K_i,\Omega_C)}=0.
\end{equation*}
Since $K_i$ converges to $K$, we from Proposition \ref{Weakly-convergence} conclude that
\[
\lambda_\phi(K_i,\eta)\to \lambda_\phi(K,\eta)\ \text{and}\ \lambda_\phi(K_i,\Omega_C)\to \lambda_\phi(K,\Omega_C)\ \text{as}\ i\to+\infty.
\]
Therefore, we define the quantity in the form of zero to zero by
\begin{equation*}
\frac{\lambda_\phi(K,\eta)}{\lambda_\phi(K,\Omega_C)}:=\lim_{i\to+\infty}\frac{\lambda_\phi(K_i,\eta)}{\lambda_\phi(K_i,\Omega_C)}
=0=\frac{\mu(\eta)}{\mu(\Omega_C)}.
\end{equation*}

For any $\omega\subset\Omega_C$, we define $\omega_i=\omega\cap\eta_i$, then $\omega_i\to\omega$ as $i\to+\infty$. Due to the continuity of the measures, we naturally define
\begin{equation*}
\frac{\lambda_\phi(K,\omega)}{\lambda_\phi(K,\Omega_C)}:=\lim_{i\to+\infty}\frac{\lambda_\phi(K,\omega_i)}{\lambda_\phi(K,\Omega_C)}
=\lim_{i\to+\infty}\frac{\mu(\omega_i)}{\mu(\Omega_C)}=:\frac{\mu(\omega)}{\mu(\Omega_C)}.
\end{equation*}
This completes the proof of this theorem.
\end{proof}

As a consequence, we obtain the corresponding result for the $L_p$ case.
\begin{corollary}
Let $\mu$ be a nonzero and locally finite Borel measure on $\Omega_C$, then there exists a $C$-full set $K$ such that
\begin{equation*}
\frac{\lambda_p(K,\cdot)}{\lambda_p(K,\Omega_C)}=\frac{\mu}{\mu(\Omega_C)}\ \ \mbox{on}\ \ \Omega_C.
\end{equation*}
\end{corollary}

\section{Spherical optimal transport and the proof of Corollary \ref{The-corollary-3}}\label{Section Spherical optimal transport}

In this section, we consider the spherical optimal transport related to the Orlicz-Gauss image problem for pseudo-cones and finish the proof of Corollary \ref{The-corollary-3}. Let $K\in ps(C)$ and $\mu$ be a nonzero Borel measure on $\Omega_C$, we call $T$ a $(K,\phi)$-transport map from $\Omega_{C^\circ}$ to $\Omega_C$ if $\mu$ is the push-forward of $\phi(\rho_K)\lambda$, i.e., $T_\sharp\,\phi(\rho_K)\lambda=\mu$. Denote
\begin{equation*}
\mathcal{T}_K^\phi=\{T:T\ \mbox{is}\ (K,\phi)\mbox{-transport map from}\ \Omega_{C^\circ}\ \mbox{to}\ \Omega_C\}.
\end{equation*}
Recall that the cost function $c:\Omega_{C^\circ}\times\Omega_C\to(0,+\infty)$ in \eqref{The-cost-function} is given by
\begin{equation*}
c(v,u)=\log\frac{1}{|\langle v,u\rangle|},\ \ \mbox{for}\ \ (v,u)\in\Omega_{C^\circ}\times\Omega_C.
\end{equation*}
Then, we prove Corollary \ref{The-corollary-3}, i.e., there exists a $C$-full set $K$ and a constant $\beta>0$ such that the total cost
\[
\int_{\Omega_{C^\circ}}c(v,T(v))\,\mathrm{d}\lambda(v),\ T\in\mathcal{T}_K^{\widetilde{\phi}},
\]
attains its maximum at the reverse radial Gauss map $\alpha^*_K$, where $\widetilde{\phi}=\beta\phi$.

\begin{proof}[Proof of Corollary \ref{The-corollary-3}]
By Theorem \ref{The-Main-Theorem}, there exists a $C$-full set $K$ and a constant $\beta>0$ such that $\beta\lambda_\phi(K,\cdot)=\mu$, i.e.,
\begin{equation*}
(\alpha^*_K)_\sharp\,\widetilde{\phi}(\rho_K)\lambda=\mu.
\end{equation*}
where $\widetilde{\phi}=\beta\phi$. By \eqref{Formula-need-OT}, we have
\[
\overline{h}_K(v)\leqslant|\langle u,v\rangle|\rho_K(u),\ \mbox{for any}\ (v,u)\in\Omega_{C^\circ}\times\Omega_C.
\]
Here equality holds if $v$ is a normal vector of $K$ at $\rho_K(u)u$. Taking the logarithm of the above inequality, we have
\[
c(v,u)\leqslant\log\rho_K(u)-\log\overline{h}_K(v).
\]

Therefore, for any $(K,\widetilde{\phi})$-transport map $T$, we have
\[
c(v,T(v))\leqslant\log\rho_K(T(v))-\log\overline{h}_K(v).
\]
This shows that
\[
\int_{\Omega_{C^\circ}}c(v,T(v))\,\mathrm{d}\lambda(v)\leqslant\int_{\Omega_{C^\circ}}\log\rho_K(T(v))\,\mathrm{d}\lambda(v)
-\int_{\Omega_{C^\circ}}\log\overline{h}_K(v)\,\mathrm{d}\lambda(v).
\]
Since $T_\sharp\,\widetilde{\phi}(\rho_K)\lambda=\mu$, for any bounded measurable function $f$ on $\Omega_C$, there holds
\[
\int_{\Omega_C}f(u)\,\mathrm{d}\mu(u)=\int_{\Omega_{C^\circ}}f(T(v))\widetilde{\phi}(\rho_K(T(v)))\,\mathrm{d}\lambda(v).
\]
We choose
\[
f(u)=\frac{\log\rho_K(u)}{\widetilde{\phi}(\rho_K(u))},
\]
then
\[
\int_{\Omega_{C^\circ}}c(v,T(v))\,\mathrm{d}\lambda(v)\leqslant\int_{\Omega_C}\frac{\log\rho_K(u)}{\widetilde{\phi}(\rho_K(u))}
\,\mathrm{d}\mu(u)-\int_{\Omega_{C^\circ}}\log\overline{h}_K(v)\,\mathrm{d}\lambda(v).
\]
Note that if $T=\alpha^*_K$, then equality holds in the above inequality. Thus, we have
\[
\int_{\Omega_{C^\circ}}c(v,T(v))\,\mathrm{d}\lambda(v)\leqslant\int_{\Omega_{C^\circ}}c(v,\alpha^*_K(v))\,\mathrm{d}\lambda(v).
\]
That is, the total cost
\[
\int_{\Omega_{C^\circ}}c(v,T(v))\,\mathrm{d}\lambda(v),\ T\in\mathcal{T}_K^{\widetilde{\phi}},
\]
attains its maximum at the reverse radial Gauss map $\alpha^*_K$.
\end{proof}

Moreover, if we consider $\phi(t)=t^p$, then we get the following result. We call $T$ a $(K,p)$-transport map from $\Omega_{C^\circ}$ to $\Omega_C$ if $\mu$ is the push-forward of $\rho^p_K\lambda$, i.e., $(T^\phi_K)_\sharp\rho^p_K\lambda=\mu$. Denote by $\mathcal{T}_K^p$ the class of all $(K,p)$-transport maps from $\Omega_{C^\circ}$ to $\Omega_C$.
\begin{corollary}
Let $\mu$ be a nonzero finite Borel measure on $\Omega_C$ and $0\neq p\in\mathbb{R}$, then there exists a $C$-full set $K$ such that the total cost
\[
\int_{\Omega_{C^\circ}}c(v,T(v))\,\mathrm{d}\lambda(v),\ T\in\mathcal{T}_K^p,
\]
attains its maximum at the reverse radial Gauss map $\alpha^*_K$.
\end{corollary}
\begin{proof}
By Corollary \ref{The-corollary-2}, there exists a $C$-full set $K$ such that $\lambda_p(K,\cdot)=\mu$, i.e., $(\alpha^*_K)_\sharp\,\rho^p_K\lambda=\mu$. Then, similar to the proof of Corollary \ref{The-corollary-3}, it is easy to check that the total cost
\[
\int_{\Omega_{C^\circ}}c(v,T(v))\,\mathrm{d}\lambda(v),\ T\in\mathcal{T}_K^p,
\]
attains its maximum at the reverse radial Gauss map $\alpha^*_K$.
\end{proof}


\vskip 5mm \noindent {\bf Author contributions:} All authors contributed equally to this work.

\end{document}